\documentclass[11pt]{article}
\makeindex
\usepackage[all]{xy}
\usepackage{graphicx}

\usepackage{amsmath}
\usepackage{amsfonts}
\usepackage{amssymb}
\usepackage{amsthm}

\newcommand{\isom}{\cong}

\newcommand{\Z}{{\bf{Z}}}

\newcommand{\Q}{{\bf{Q}}}

\newcommand{\R}{{\bf{R}}}
\newcommand{\C}{{\bf{C}}}
\newcommand{\F}{{\bf{F}}}
\newcommand{\T}{{\bf{T}}}

\newcommand{\m}{{\mathfrak{m}}}
\newcommand{\PP}{\bf{P}}

\newcommand{\Mid}{|} 

\newcommand{\miD}{|}

\newcommand{\divs}{\!\mid\!}
\newcommand{\ndiv}{\!\nmid\!}
\newcommand{\tensor}{\otimes}

\newcommand{\ra}{{\rightarrow}}
\newcommand{\lra}{{\longrightarrow}}

 1
\DeclareFontEncoding{OT2}{}{} 
  \newcommand{\textcyr}[1]{%
    {\fontencoding{OT2}\fontfamily{wncyr}\fontseries{m}\fontshape{n}%
     \selectfont #1}}
\newcommand{\Sha}{{\mbox{\textcyr{Sh}}}} 
\newcommand{\ansha}{{\Mid \Sha(A_f) \miD_{\rm an}}}

\newcommand{\If}{{I_{\scriptscriptstyle{f}}}}

\newcommand{\ED}{{E_{\scriptscriptstyle{-D}}}}
\newcommand{\omegaED}{{\omega_{\scriptscriptstyle{E_{-D}}}}}
\newcommand{\omegaE}{{\omega_{\scriptscriptstyle{E}}}}

\newcommand{\cF}{{\mathfrak{c}}(A_f/F)}
\newcommand{\cK}{{\mathfrak{c}}(A_f/K)}
\newcommand{\qq}{{\mathfrak{q}}}
\newcommand{\mm}{{\mathfrak{m}}}

\newcommand{\ord}{{\rm ord}}
\newcommand{\LAf}{{L(A_f,1)}}

\newcommand{\OAf}{{\Omega(A_f)}}
\newcommand{\OAg}{{\Omega(A_g)}}
\newcommand{\OAfm}{{\Omega^-_{\scriptscriptstyle{A_f}}}}
\newcommand{\cAf}{{c_{\scriptscriptstyle{A_f}}}}
\newcommand{\cA}{{c_{\scriptscriptstyle{A}}}}
\newcommand{\cE}{{c_{\scriptscriptstyle{E}}}}
\newcommand{\OAfmD}{{\OAfm / (-D)^{d/2}}}
\newcommand{\OEm}{{\Omega^-_{\scriptscriptstyle{E}}}}
\newcommand{\OAfp}{{\Omega^+_{\scriptscriptstyle{A_f}}}}
\newcommand{\sD}{{(\sqrt{-D})}}
\newcommand{\NerA}{{\mathcal{A}}}

\newcommand{\sP}{{\mathcal{P}}}
\newcommand{\OO}{{\mathcal{O}}}
\newcommand{\OD}{{\mathcal{O}_{-D}}}
\newcommand{\Pz}{{\mathcal{P}^0}}
\newcommand{\Pzm}{{\mathcal{P}^0_\mm}}
\newcommand{\Hpm}{{H^+_\mm}}
\newcommand{\Hmm}{{H^-_\mm}}
\newcommand{\Tm}{{\T_\mm}}
\newcommand{\tent}{{\tensor_{\T}}}
\newcommand{\tenth}{{\tensor_{\T[1/2]}}}

\newcommand{\tentq}{{\tensor_{\T_\Q}}}
\newcommand{\tentm}{{\tensor_{\T_\mm}}}

\newcommand{\eD}{{e_{\scriptscriptstyle{D}}}}
\newcommand{\chiD}{{\chi^0_{\scriptscriptstyle{D}}}}

\newcommand{\epD}{{\epsilon_{\scriptscriptstyle{D}}}}

\newcommand{\comment}[1]{}

\newtheorem{lem}{Lemma}[section]
\newtheorem{cor}[lem]{Corollary}
\newtheorem{prop}[lem]{Proposition}
\newtheorem{conj}[lem]{Conjecture}
\newtheorem{thm}[lem]{Theorem}

\theoremstyle{definition}

\newtheorem{rmk}[lem]{Remark}

\begin{document}
\markboth{Amod Agashe}
{Squareness in the special $L$-value}

\title{Squareness in the special $L$-value \\ or\\
Squareness in the special $L$-value and special $L$-values
of twists
\footnote{
I would like to modify the title and include the words
``and special $L$-values of twists''. I have included the earlier title
in case it was used as an identifier for the paper.}
}

\author{Amod Agashe\footnote{
The author was partially supported by National Science Foundation Grant No. 0603668.}}

\maketitle





\begin{abstract}
Let $N$ be a prime and 
let~$A$ be a quotient of~$J_0(N)$ over~$\Q$ associated to a newform
such that the special $L$-value of~$A$ (at $s=1$) is non-zero.
Suppose that the algebraic part of the special $L$-value of~$A$ is divisible
by an odd prime~$q$ such that $q$
does not divide the numerator of~$\frac{N-1}{12}$.
Then the Birch and Swinnerton-Dyer conjecture predicts
that the $q$-adic valuations of 
the algebraic part of the special $L$-value of~$A$
and of the 
order of the Shafarevich-Tate group are both positive even numbers. 
Under a certain
\mbox{$\bmod\ q$} non-vanishing hypothesis on special $L$-values of twists
of~$A$, we show 
that the \mbox{$q$-adic} valuations of 
the algebraic part of the special $L$-value of~$A$ and of the Birch and 
Swinnerton-Dyer conjectural
order of the Shafarevich-Tate group of~$A$ are both positive even 
numbers.
We also give a formula for the algebraic part of the special
$L$-value of~$A$ over quadratic imaginary fields~$K$
in terms of the free abelian group on isomorphism classes of supersingular
elliptic curves in characteristic~$N$ (equivalently, over conjugacy classes
of maximal orders in the definite quaternion algebra over~$\Q$ ramified at~$N$
and~$\infty$)
which shows that this algebraic part is a 
perfect square up to powers of the prime two and of
primes dividing the discriminant
of~$K$.
Finally, for an optimal elliptic curve~$E$, 
we give a formula for the special $L$-value of the twist~$\ED$ of~$E$
by a negative fundamental discriminant~$-D$, which shows
that this special $L$-value is an integer up to a power of~$2$, 
under some hypotheses.
In view of the second part of the Birch and Swinnerton-Dyer conjecture,
this leads us to the surprising conjecture that the square of
the order of the torsion subgroup of~$\ED$ divides the product
of the order of the Shafarevich-Tate group of~$\ED$ and the orders
of the arithmetic component groups of~$\ED$, under certain mild hypotheses.
\end{abstract}

\section{Introduction} \label{sec:intro}


Let $A$ be an abelian variety over a number field~$F$, and 
let $L(A/F,s)$ denote the associated $L$-function,
which we assume is defined over all of~$\C$ (this will be true
in the cases we are interested in).
Let $\Omega(A/F)$ denote the quantity~$C_{A,\infty}$ 
in~\cite[\S~III.5]{lang:nt3};
it is the ``archimedian volume'' of~$A$ over embeddings of~$F$ in~$\R$ and~$\C$
(e.g., if $F = \Q$, then it is the volume of~$A(\R)$ computed using
a generator for the highest exterior power of the group of
invariant differentials on the N\'eron model of~$A$;
the only other case we shall need is when $F$ is a quadratic imaginary
field,  which is discussed at the beginning of Section~\ref{sec:overK}).
Let $M_{\rm fin}$ denote the set of finite places of~$F$.
Let $\NerA$ denote the N\'eron model of~$A$ over the ring of integers
of~$F$ and let $\NerA^0$ denote 
the largest open subgroup scheme
of~$\NerA$ in which all the fibers are connected.
If $v \in M_{\rm fin}$, 
then let $\F_v$ denote the associated residue class field and let 
$c_v(A/F) = [\NerA_{\F_v}(\F_v): \NerA^0_{\F_v}(\F_v)]$, the orders
of the arithmetic component groups.
Let $\Sha(A/F)$ denote
the Shafarevich-Tate group of~$A$ over~$F$.
If $F = \Q$, then we will often drop the symbol ``$/F$''
in the notation
(thus $\Sha(A/\Q)$ will be denoted $\Sha(A)$, etc.).
If $B$ is an abelian variety over~$F$, then we denote by~$B^{\vee}$
the dual abelian variety of~$B$, and by~$B(F)_{\rm tor}$ the torsion
subgroup of~$B(F)$. Suppose that 
$L(A/F,1) \neq 0$. 
Then
the second part of the Birch and Swinnerton-Dyer conjecture says the following
(see~\cite[\S~III.5]{lang:nt3}):

\begin{conj}[Birch and Swinnerton-Dyer] \label{bsd}
\begin{eqnarray} \label{bsdform1}
\frac{L(A/F,1)}{\Omega(A/F)} =
\frac {\Mid \Sha(A/F) \miD \cdot \prod_{\scriptscriptstyle{v \in M_{\rm fin}}}  
c_v(A/F)}
      { \Mid A(F)_{\rm tor} \miD \cdot \Mid A^{\vee}(F)_{\rm tor} \miD }.
\end{eqnarray}
\end{conj}
We denote by~$\Mid \Sha(A/F) \miD_{\rm an}$ the value of
$\Mid \Sha(A/F) \miD$ 
predicted by the conjecture above, and call it the
{\em analytic order} of~$\Sha(A/F)$. Thus
$$\Mid \Sha(A/F) \miD_{\rm an}
= \frac{L(A/F,1)}{\Omega(A/F)} \cdot 
\frac{ \Mid A(F)_{\rm tor} \miD \cdot \Mid A^{\vee}(F)_{\rm tor} \miD }
{\prod_{\scriptscriptstyle{v \in M_{\rm fin}}}  c_v(A/F)}\ \ .
$$
We also call the ratio $\frac{L(A/F,1)}{\Omega(A/F)}$, 
the {\em algebraic part} of the special
$L$-value of~$A_f$ over~$F$; in the contexts where we shall use this,
it is known that the ratio is a rational number (and hence an
algebraic number).

If $N$ is a positive integer, then
let $X_0(N)$ denote the modular curve over~$\Q$ associated to~$\Gamma_0(N)$,
and let $J_0(N)$ be its Jacobian.
Let $\T$ denote the subring of endomorphisms of~$J_0(N)$
generated by the Hecke operators (usually denoted~$T_\ell$
for $\ell \ndiv N$ and $U_p$ for $p\divs N$). 
If $f$ is a newform of weight~$2$ on~$\Gamma_0(N)$, then 
let $I_f = {\rm Ann}_{\T} f$ and 
let $A_f$ denote the quotient abelian variety $J_0(N)/ I_f J_0(N)$
over~$\Q$. 
We also denote by~$L(f,s)$ the $L$-function associated to~$f$ and
by~$L(A_f,s)$ the $L$-function associated to~$A_f$.
It is known that $\frac{L(A_f,1)}{\OAf}$ is a rational number.

Now fix 
a newform~$f$ of weight~$2$ on $\Gamma_0(N)$ such
that $L(A_f,1) \neq 0$. Then 
by~\cite{kollog:finiteness}, 
$A_f(\Q)$ has rank zero
and $\Sha(A_f)$ is finite.
Thus
the second part of the Birch and Swinnerton-Dyer conjecture becomes:
\comment{
, where $\Sha(A_f)$ denotes
the Shafarevich-Tate group of~$A_f$.
Let $\NerA$ denote the N\'eron model of~$A_f$ over~$\Z$
and let $\NerA^0$ denote the largest open subgroup scheme
of~$\NerA$ in which all the fibers are connected.
Let $d = \dim A_f$, and let $\omega$ be a generator of
the $d$-th exterior power of the group of invariant
differentials on~$\NerA$.
Let $\OAf$ denote the volume of~$A_f(\R)$ with
respect to the measure given by~$\omega$.
If $p$ is a prime number, then let 
$c_p(A) = [\NerA_{\F_p}(\F_p): \NerA^0_{\F_p}(\F_p)]$.
If $A$ is an abelian variety, then we denote by~$A^{\vee}$
the dual abelian variety of~$A$.
}

\begin{conj}[Birch and Swinnerton-Dyer] \label{bsd2}
\begin{eqnarray} \label{bsdform}
\frac{L(A_f,1)}{\OAf} =
\frac {\Mid \Sha(A_f) \miD \cdot \prod_{p | N} c_p(A_f)}
      { \Mid A_f(\Q) \miD \cdot \Mid A_f^{\vee}(\Q) \miD },
\end{eqnarray}
\end{conj}

Recall that an integer is said to be a 
{\em fundamental discriminant} if it is the discriminant
of a quadratic field. 
The results of this paper concern the algebraic parts of
the special $L$-values of~$A_f$ over~$\Q$, of~$A_f$ over
quadratic imaginary fields, and of twists of~$A_f$ by
negative fundamental discriminants (over~$\Q$). 
In Section~\ref{sec:twists}, when $A_f$ is an elliptic curve,
we give 
a formula for the special $L$-value of the twist of~$A_f$
by a negative fundamental discriminant, which shows
that this special $L$-value is an integer, under some hypotheses.
This leads us to 
the surprising conjecture that for such twists, the {\em square} of
the order of the torsion subgroup divides the product
of the order of the Shafarevich-Tate group and the orders
of the arithmetic component groups, under certain mild hypotheses.
In Section~\ref{sec:overQ}, 
under a certain
\mbox{$\bmod\ q$} non-vanishing hypothesis on special $L$-values of twists
of~$A_f$, we show that when $N$ is prime,
for certain odd primes~$q$ that divide the algebraic part of the 
special $L$-value of~$A_f$ over~$\Q$, the \mbox{$q$-adic} valuations of 
the algebraic part of the special $L$-value of~$A_f$ and of the Birch and 
Swinnerton-Dyer conjectural
order of the Shafarevich-Tate group of~$A_f$ are both positive even 
numbers, in conformity with what the second part of the Birch and
Swinnerton-Dyer conjecture predicts.
In Section~\ref{sec:overK}, for $N$ prime,
we give a formula for the algebraic part of the special
$L$-value of~$A_f$ over quadratic imaginary fields~$K$
in terms of the free abelian group on isomorphism classes of supersingular
elliptic curves in characteristic~$N$ (equivalently over conjugacy classes
of maximal orders in the definite quaternion algebra over~$\Q$ ramified at~$N$
and~$\infty$)
which shows that this algebraic part is a 
perfect square away from the prime two and 
the primes dividing the discriminant
of~$K$. In Section~\ref{sec:proofs}, we give the proofs of two 
theorems mentioned in Sections~\ref{sec:overQ}  and~\ref{sec:overK}. 
Finally, in Section~\ref{sec:app},
we we give a formula for the determinant of 
the ``complex period matrix'' of an abelian variety, which
is needed in the proof of the main theorem of Section~\ref{sec:overK}.
All the sections except Section~\ref{sec:proofs}
may be read independently of each other, although there is some
cross referencing.

We now introduce some notation that will be used in various sections
of this article.
If $\langle\ ,\ \rangle : M \times M' \ra \C$, 
is a pairing between two $\Z$-modules~$M$ and~$M'$, each
of the same rank~$m$, and 
$\{\alpha_1, \ldots, \alpha_m \}$
and $\{ \beta_1, \ldots, \beta_m\}$ are 
bases of~$M$ and~$M'$ (respectively),
then by~${\rm disc}(M \times M' \ra \C)$, we mean
$\det (\langle \alpha_i, \beta_j\rangle )$. 
Up to a sign, ${\rm disc}(M \times M' \ra \C)$ is independent 
of the choices of bases made in its definition,
and in the rest of this paper, ${\rm disc}(M \times M' \ra \C)$
will be well defined only up to a sign
(this ambiguity will not matter for our main results).
We have a pairing
\begin{eqnarray} \label{eqn:pairing}
H_1(X_0(N),\Z) \tensor {\C} \times  S_2(\Gamma_0(N), {\C}) \ra {\C}
\end{eqnarray}
given by~$(\gamma,g)\mapsto \langle \gamma,g\rangle  = \int_\gamma 2 \pi i g(z) dz$
and extended $\C$-linearly.
At various points in this article, we will consider pairings 
between two $\Z$-modules; unless otherwise stated,
each such pairing is obtained in a natural
way from~(\ref{eqn:pairing}).
We have an involution induced by complex conjugation on
$H_1(A_f,\Z)$. Let $H_1(A_f,\Z)^+$ and $H_1(A_f,\Z)^-$ denote
the subgroups
of elements of~$H_1(A_f,\Z)$ on which the involution acts as~$1$
and~$-1$ respectively.
Let $S_f = S_2(\Gamma_0(N),\Z)[I_f]$\label{sym:sf}, 
let $\OAfp = {\rm disc}(H_1(A_f,\Z)^+ \times S_f \ra {\C})$, and 
let $\OAfm = {\rm disc}(H_1(A_f,\Z)^- \times S_f \ra {\C})$.
In each section below, we will continue to use the notation 
introduced in this section, unless mentioned otherwise.

\section{Special $L$-values of twists of elliptic curves} \label{sec:twists}

In this section, we give a formula for the special $L$-value
of the twist of an optimal elliptic curve by 
a negative fundamental discriminant,
which shows that this special $L$-value is an integer up to a power of~$2$,
under certain hypotheses. This has some surprising implications from 
the point of view of the Birch and Swinnerton-Dyer conjecture, as we shall 
discuss.

%

We now recall some definitions for an elliptic curve~$A$ defined
over~$\Q$.
If $d$ is a square free integer, 
then $A_d$ denotes the
twist of~$A$ by~$d$. Thus if $y^2 = x^3 + ax +b$ with $a, b \in \Q$
is a Weierstrass equation for~$A$, then 
$y^2 = x^3 + d^2 a x + d^3 b$ is a Weierstrass equation for~$A_d$.
If $-D$
is a negative fundamental discriminant, we shall often consider
the following hypothesis on~$(A,-D)$:\\
(**) $-D$ is coprime to the 
discriminant of some
Weierstrass 
equation $y^2 = x^3 + Ax + B$ for~$E$ with $A, B \in \Z$.

Note that for every elliptic curve over~$\Q$, there is a Weierstrass 
equation $y^2 = x^3 + Ax + B$ with $A, B \in \Z$.
If $A$ is an elliptic curve over~$\Q$, then let $\omega_A$ denote
an invariant differential on a global minimal Weierstrass model of~$A$,
which is unique upto sign. Now assume that $A$ is an optimal elliptic curve,
i.e., it is~$A_f$ for some newform~$f$ of weight~$2$ on~$\Gamma_0(N)$
for some~$N$. Let $\pi: X_0(N) \ra A$ denote the associated parametrization.
Then the space of pullbacks of differentials on~$A$ to~$X_0(N)$
is spanned by the differential~$2 \pi i f(z) dz$; let 
$\omega_f$ be the differential on~$A$ whose pullback
is precisely~$2 \pi i f(z) dz$.  
Then $\omega_A = c \omega_f$ for some rational number~$\cA$,
which is called the Manin constant of~$A$.

\begin{lem} \label{lem:twist}
Let $E$ be an optimal elliptic curve over~$\Q$ and 
let $-D$ be a negative fundamental discriminant such that 
$(E, -D)$ 
satisfies hypothesis~(**).
Then up to a sign,
$$\Omega(\ED) = \cE \cdot c_\infty(\ED) \cdot   \OEm  / \sqrt{-D}\ ,$$
where $c_\infty(\ED)$ is the number of connected
components of~$\ED(\R)$.
\end{lem}
\begin{proof}
By hypothesis~(**), 
there is a 
Weierstrass 
equation $y^2 = x^3 + Ax + B$ for~$E$ with $A, B \in \Z$, such that 
$-D$ is coprime to the 
discriminant of this equation.
Denote this equation  by~$(a)$.
If $(x)$ denotes a Weierstrass equation for an elliptic curve, then
we denote the associated discriminant by~$\Delta(x)$ and
the associated invariant differential by~$\omega(x)$.
Replacing $x$ by~$x/\sD^2$ and $y$ by~$y/\sD^3$,
we get the Weierstrass equation $y^2 = x^3 + D^2 A x - D^3 B$
for~$\ED$ (in fact, this transformation gives an isomorphism
of~$E$ and~$\ED$ over~$\Q(\sqrt{-D})$); denote this equation by~$(b)$.
Then by~\cite[Table~III.1.2]{silverman:aec},
$\Delta(b) = D^6 \Delta(a)$ and 
\begin{eqnarray} \label{eqn:omega}
\omega(b) = \omega(a) /\sD.
\end{eqnarray}
Since $D$ is squarefree and coprime to~$\Delta(a)$, if $p$ is a prime
that divides~$D$,  then $\ord_p(\Delta(a)) = 0 <12$,
and $\ord_p(\Delta(b)) = \ord_p(D^6 \Delta(a)) = 6 <12$. Thus
by~\cite[Rmk.~VII.1.1]{silverman:aec},
equations~$(a)$ and~$(b)$ are both minimal at the primes dividing~$D$.
Also, if $p$ is a prime that does not divide~$D$,
then the coefficients of~$(a)$ and~$(b)$ have the same order at~$p$.
Thus, following 
the proof of Prop.~VIII.8.2 in~\cite{silverman:aec}, 
there is a transformation
$x = u^2 x' + r$, $y = u^3 y' + u^2 s x' + t$ for some integers
$u, r, s$ and~$t$, with $u \neq 0$, which
converts both equations~$(a)$ and~$(b)$ to equations that are minimal at 
all primes. Denote these equations by~$(c)$ and~$(d)$ respectively;
these are then global minimal Weierstrass equations for~$E$ and~$\ED$
respectively. Hence $\omegaE = \omega(c)$ and $\omegaED = \omega(d)$.
By~\cite[p.~49]{silverman:aec},
$\omega(c) = u \omega(a)$
and $\omega(d) = u \omega(b)$.
Using~(\ref{eqn:omega}),
$\omega(d) = u \omega(b) = u  \omega(a) / \sD
= \omega(c) /\sD$. 

Also, equation~$(b)$ was obtained from equation~$(a)$ by
replacing $x$ by~$x/\sD^2$ and $y$ by~$y/\sD^3$. 
Thus if $(x,y)$ is a point on~$(b)$, then
the corresponding point on~$(a)$ is given by~$(x/\sD^2, y/\sD^3)$.
Since the transformation used to go from~$(b)$ to~$(d)$ was the
same as the one used to go from~$(a)$ to~$(c)$, we see
that again, 
if $(x,y)$ is a point on~$(d)$, then
the corresponding point on~$(c)$ is given by~$(x/\sD^2, y/\sD^3)$.
Denote this map from points on~$(b)$ to points on~$(a)$ by~$T$ and 
let $\sigma$ denote complex conjugation.
Then if $P = (x, y)$ is a point on~$(d)$ that 
is fixed by complex conjugation, i.e. $\sigma(x,y) = (x,y)$,
then $\sigma(T(P)) = \sigma(x/\sD^2, y/\sD^3) = (x/\sD^2, - y/\sD^3)
= - T(P)$. From this we see that if $\gamma \in H_1(\ED, \Z)$ is a generator, 
then $T(\gamma) \in H_1(E, \Z)^-$. It is easy to see that $T$ is
invertible, and so $T(\gamma)$ is a generator of~$H_1(E, \Z)^-$. 

Thus $\OEm = \int_{T(\gamma)} \omega_f$ up to a sign, and 
using the change of variables given by the transformation~$T$,
we see that 
$\int_\gamma \omegaED = \int_{T(\gamma)} \omegaE / \sD$.
Also, recall that
$\omegaE = \cE \omega_f$.
From the discussion above, 
we see that up to a sign,
$\Omega(\ED)
=  c_\infty(\ED)  \cdot  \int_\gamma \omegaED 
=   c_\infty(\ED)  \cdot \int_{T(\gamma)} \omegaE / \sD 
=  c_\infty(\ED)  \cdot \cE \cdot \int_{T(\gamma)} \omega_f / \sD 
=  c_\infty(\ED)  \cdot \cE \cdot \OEm / \sqrt{-D}$, as was to be shown.

\end{proof}

Let $N$ be a positive integer and let $f$ be a newform of 
weight~$2$ on~$\Gamma_0(N)$.
Let $-D$ be a negative fundamental discriminant that is
coprime to~$N$ and 
let $\epD = 
(\frac{-D}{\cdot})$ denote the quadratic character associated to~$-D$.
If $f(q) = \sum_{n > 0} a_n q^n$ is the Fourier expansion of~$f$,
then the twist of~$f$ by~$\epD$ is the modular form whose Fourier expansion is
$(f \tensor \epD)(q) = \sum_{n > 0} \epD(n) a_n q^n$. It is in fact
a newform in~$S_2(ND^2, \epsilon_{\scriptscriptstyle{D}}^2)$ (considering that $D$ is coprime to~$N$; see, e.g., p.~221 and p.~228 of~\cite{atkin-li}
and the references in loc. cit.).
Just as we associated an abelian variety~$A_f$ to~$f$, one can associate
to~$f \tensor \epD$ an abelian variety quotient~$A_{f \tensor \epD}$ of~$J_1(ND^2)$,
and 
moreover, if $f_1, \ldots, f_d$ are the Galois conjugates of~$f$, then
$L(A_{f \tensor \epD},1) = \prod_i L(f_i \tensor \epD,1)$
(see, e.g. p.~89 and p.~95 of~\cite{rohrlich-modular}).

\begin{prop} \label{prop:twist}
Suppose $f$ has integer Fourier coefficients, and
let $E$ denote the associated optimal elliptic curve quotient of~$J_0(N)$ 
over~$\Q$.
Suppose that 
$(E, -D)$ 
satisfies hypothesis~(**) mentioned at the beginning of this section. Then
$$\frac{L(\ED,1)}{\Omega(\ED)}
= \frac{L(A_{f \tensor \epD},1)}{\cE \cdot
c_\infty(\ED) \cdot \OAfm / \sqrt{-D}}\ ,$$ 
where recall that 
$c_\infty(\ED)$ is the number of connected
components of~$\ED(\R)$, 
$\cE$ is the Manin constant of~$E$,
and $\OAfm$ is as defined at the end of Section~\ref{sec:intro}.  
In particular, if $N$ is square free or if $\cE = 1$
(as is conjectured), then
$$\frac{L(A_{f \tensor \epD},1)}{\OAfm / \sqrt{-D}}
= \frac{L(\ED,1)}{\Omega(\ED)}\ ,$$ 
up to a power of~$2$.
\end{prop}
\begin{proof}
The first statement  follows from Lemma~\ref{lem:twist} above,
considering that 
$L(\ED,1) = L(f\tensor \epD,1) =  L(A_{f \tensor \epD},1)$. The second statement
follows from the first,
considering that $c_\infty(\ED)$ is a power of~$2$, and 
if $N$ is squarefree, then
$\cE$ is a power of~$2$ as well
(by~\cite[Cor.~4.1]{mazur:rational}). 
\end{proof}

The modular symbol
$\sum_{\scriptscriptstyle{b \bmod D}} \epsilon_{\scriptscriptstyle D}(b) \{ - \frac{b}{D},\infty \}$
is an element of~$H_1(X_0(N), \Z)^-$ by~\cite[Lemma~5.2]{marusia:module} 
(see also~\cite[\S9.8--9.9]{manin:cyclo}),
and will be denoted by~$\eD$.
Let $\pi$ denote the quotient map $J_0(N) \ra A_f$, and 
let $\pi_*$ denote the induced map $H_1(J_0(N)(\C), \Q) \ra H_1(A_f(\C), \Q)$.
Let $d = \dim A_f$.

\begin{prop}\label{prop:ftensore}
Assume that 
$L(A_{f \tensor \epD},1) \neq 0$. 
Then up to a power of~$2$,
$$\frac{L(A_{f \tensor \epD},1)}{\OAfmD}
= | \pi_* (H_1(X_0(N), \Z)^-): \pi_*(\T \eD) | \ .$$
\end{prop}

\begin{proof}
The proof is very similar to the proof of Theorem~2.1 in~\cite{agmer}.
The main thing to note is that
if $f_1, \ldots, f_d$ are the Galois conjugates of~$f$, then
for $i = 1, \ldots, d$, we have
$L(f_i \tensor \epD,1) = \frac{\langle \eD,f_i\rangle}{\sqrt{-D}}$
(see, e.g., \cite[p.~254]{marusia:module} or~\cite[Thm~9.9]{manin:cyclo}),
and so
$L(A_{f \tensor \epD},1) = \prod_i L(f_i \tensor \epD,1) 
= \prod_i \frac{\langle \eD,f_i\rangle}{\sqrt{-D}}$.
Also, up to a power of~$2$, $\pi_*(H_1(X_0(N), \Z)^-) = H_1(A_f, \Z)^-$.
Hence, up to a power of~$2$, 
\begin{eqnarray*}
\frac{L(A_{f \tensor \epD},1)}{\OAfm / (-D)^{d/2} }
& = &\frac {\prod_i \langle \eD,f_i\rangle }
{{\rm disc}(\pi_*(H_1(X_0(N), \Z)^-) \times S_f \ra {\C})} \\
& = &\frac {\prod_i \langle \eD,f_i\rangle }
{{\rm disc}(\pi_*(\T \eD) \times S_f \ra {\C})} \cdot
| \pi_*(H_1(X_0(N), \Z)^-) : \pi_*({\T}\eD)|.
\end{eqnarray*}

One can see in a manner similar to the proof of formula~(6)
in the proof of Theorem~2.1 in~\cite{agmer}
that the first factor above is~$1$, as we explain next.
The proposition then follows from the claim in the previous sentence.

There is a perfect pairing 
$$
\T \times S_2(\Gamma_0(N),\Z) \ra \Z
$$
which associates to~$(T,f)$ the first Fourier coefficient~$a_1(f \divs T)$ of 
the modular form~$f \divs T$ (see \cite[(2.2)]{ribet:modp});
this induces a pairing
\begin{eqnarray}\label{perfpair}
\psi: \T/I_f \times S_f \ra \Z,
\end{eqnarray}
which is also perfect
(e.g., see~\cite[Lemma~2.2]{agmer}).

\noindent {\em Claim:} The map $\T \ra \T \eD$ given by $t \mapsto t\eD$
induces an isomorphism $\T/I_f \stackrel{\cong}{\lra}\T \eD/ \If \eD$. 

\noindent {\em Proof}:
It is clear that the map $\T \ra \T \eD/ \If \eD$ given by $t \mapsto t\eD$
is surjective. All we have to show is that
the kernel of this map is~$I_f$.
It is clear that the kernel contains~$I_f$. Conversely, 
if $t$ is in the kernel, then $t\eD \in \If \eD$; let
$i \in I_f$ be such that $t\eD = i\eD$. Then $(t-i)\eD=0$,
and thus $\int_{(t-i)\eD} \omega_f = 0$, i.e., 
$\int_\eD \omega_{(t-i)f} =0$. If the eigenvalue
of~$f$ under~$(t-i)$ is~$\lambda$, then this means
$\lambda \cdot L(f \tensor \epD,1) = 0$, i.e., $\lambda = 0$
(since $L(f \tensor \epD,1) \neq 0$, considering
that $L(A_{f \tensor \epD},1) \neq 0$). Thus
$(t-i) \in I_f$, i.e., $t \in I_f$. This proves the claim.


We continue the proof of the theorem.
In what follows,
$i, j, k, {\rm and\ } \ell$ are indices running from~$1$ to~$d$.
Let $\{g_k\}$ be a~$\Z$-basis of~$S_f$ and let $\{t_j\}$ be the 
corresponding dual basis
of~${\T}/I_f$ under the perfect pairing~$\psi$ in~(\ref{perfpair}) above. 
Then by the claim above, $\{t_j \eD\}$ is a $\Z$-basis 
for~${\T} \eD / \If \eD$.
Now $g_k = \sum_k a_{ki} f_i$ for some $\{a_{ki} \in \C\}$. Let
$A$ be the matrix having $(k,i)$-th entry~$a_{ki}$,
and let $(a^{-1})_{i\ell}$
denote the $(i,\ell)$-th element of the inverse of~$A$. Then 

\begin{tabbing}
\= ${\rm disc} ({\T} \eD / \If \eD \times S_f \ra {\C})$ \\
\> $ = \det \{\langle t_j \eD,g_k\rangle \}= \det \{\langle \eD,g_k \divs t_j\rangle \} 
   = \det \{\langle \eD,(\sum_i a_{ki} f_i)\divs t_j\rangle \}$\\ 
\> $= \det \{\langle \eD,\sum_i a_{ki} a_1(f_i \divs t_j) f_i\rangle \}$ 
   \ \ \ \ (since $f_i$'s are eigenvectors) \\
\> $= \det \{\langle \eD,\sum_i a_{ki} \sum_\ell (a^{-1})_{i\ell} a_1(g_\ell \divs t_j) f_i\rangle \}$
   \ \ (using $f_i= \sum_\ell (a^{-1})_{i\ell} g_\ell$)\\
\> $= \det \{\langle \eD,\sum_i a_{ki} (a^{-1})_{ij} f_i\rangle \}$
   \ \ \ \ (using $a_1(g_\ell \divs t_j) = \delta_{\ell j}$)\\
\> $= \det \{\sum_i a_{ki} (a^{-1})_{ij} \langle \eD,f_i\rangle \} 
= \det \{\sum_i a_{ki} \langle \eD,f_i\rangle  (a^{-1})_{ij} \}$\\ 
\> $= \det(A \Delta A^{-1})$
  \ \ \ \  (where $\Delta = {\rm diag}(\langle \eD,f_i\rangle )$)\\
\> $= \det (\Delta) = \prod_i \langle \eD,f_i\rangle .$
\end{tabbing}

This shows what we wanted and finishes the proof of the proposition.
\end{proof}

\begin{cor} \label{cor:twist}
Let $E$ be an optimal elliptic curve over~$\Q$ and 
let $-D$ be a negative fundamental discriminant such that 
$(E, -D)$ 
satisfies hypothesis~(**) mentioned at the beginning of this section. 
Assume either that the Manin constant of~$E$ is one
(as conjectured) or that $N$ is squarefree. Then 
$$\frac{L(\ED,1)}{\Omega(\ED)} \in \Z[1/2]. $$
\end{cor}
\begin{proof}
This follows from Propositions~\ref{prop:twist} and~\ref{prop:ftensore}.
\end{proof}

In view of the Birch and Swinnerton-Dyer conjecture 
(Conjecture~\ref{bsd2} above) and the conjecture that the
Manin constant is one, the corollary above suggests
the following conjecture:

\begin{conj}
Let $E$ be an optimal elliptic curve over~$\Q$ of conductor~$N$ and
let $-D$ be a negative fundamental discriminant such that 
$(E, -D)$ 
satisfies hypothesis~(**) mentioned at the beginning of this section. 
Recall that $\ED$ denotes the twist of~$E$ by~$-D$. Suppose 
$L(\ED,1) \neq 0$.
Then 
$\Mid \ED(\Q) \miD^2$ divides
${\Mid \Sha(\ED) \miD \cdot   \prod_{p | N} c_p(\ED)}$,
up to a power of~$2$.
\end{conj}

Using the mathematical software sage, with its inbuilt Cremona's database 
for all elliptic curves of conductor up to~$130000$,
we verified the conjecture above for all triples $(N, E, D)$ such 
that $N$ and~$D$ are positive integers with $N D^2 < 130000$, 
and $E$ is an optimal elliptic curve
of conductor~$N$. In fact, we found that even if replace
the hypothesis~(**)
with the potentially weaker hypothesis that 
${\rm gcd}(N,D) = 1$, the conclusion of the conjecture above
was true in all examples, even at the prime~$2$ (i.e., not just
up to a power of~$2$).
We also found that in all these examples, 
the odd part of 
$\Mid \ED(\Q) \miD^2$ divides
$\prod_{p | N} c_p(\ED)$, and that if $-D \neq -3$, then
$\Mid \ED(\Q) \miD$ is a power of~$2$.
Table~\ref{table1} below shows some interesting examples.
The example of $E = 105a1$ shows that
$\Mid \ED(\Q) \miD^2$ does not divide
$\prod_{p | N} c_p(\ED)$ in general
(but it does divide
$\Mid \Sha(\ED) \miD \cdot   \prod_{p | N} c_p(\ED)$).
Also, if $-D = -3$, 
it is not true that $\Mid \ED(\Q) \miD$ is a power of~$2$,
as the example of $E = 14a1$  shows.
If we relax the assumption that ${\rm gcd}(N,D) = 1$, then
it is no longer true that 
$\Mid \ED(\Q) \miD^2$ divides
${\Mid \Sha(\ED) \miD \cdot   \prod_{p | N} c_p(\ED)}$, 
as the examples $E = 21a1$ and $E = 27a1$ show.

\begin{table}\caption{\label{table1}}
\begin{center}
\begin{tabular}{|l|l|l|l|l|}\hline
$E$ &  $-D$ & $|\ED(\Q)_{\rm tor}|$ & $\prod_p c_p(\ED)$ & $|\Sha(\ED)_{\rm an}|$
\\\hline
14a1 & -3 &  6 & 36 & 1  \\\hline
21a1 & -7 & 4 & 8 & 1 \\\hline
27a1 & -3 & 3 & 1 & 1 \\\hline
105a1 & -11 & 2 & 2 & 4\\\hline
\end{tabular}
\end{center}
\end{table}

\section{Special L-values over~$\Q$} \label{sec:overQ}
We assume in this section  that $N$ is prime.
Let $f$ be a newform of weight~$2$ on~$\Gamma_0(N)$, and as before
let $A_f$ denote the associated newform quotient of~$J_0(N)$ 
over~$\Q$.
Let
$q$ be an odd prime that does not divide the numerator of~$\frac{N-1}{12}$
but divides $\frac{\LAf}{\OAf}$. 
Note that the denominator of~$\frac{\LAf}{\OAf}$ divides 
the numerator of~$\frac{N-1}{12}$ (e.g., by~\cite[Prop.~4.6]{agst:bsd} and
the fact that the order of the cusp~$(0) - (\infty) \in J_0(N)(\C)$ 
is the numerator of~$\frac{N-1}{12}$ when $N$ is prime), and so 
it makes sense to talk about whether $q$ divides $\frac{\LAf}{\OAf}$ or not.
In this section, we show that 
under a certain
\mbox{$\bmod\ q$} non-vanishing hypothesis on special $L$-values of twists
of~$A_f$, the \mbox{$q$-adic} valuations of 
the algebraic part of the special $L$-value of~$A_f$ and of the Birch and 
Swinnerton-Dyer conjectural
order of the Shafarevich-Tate group of~$A_f$ are both positive even 
numbers, in conformity with what the second part of the Birch and
Swinnerton-Dyer conjecture predicts.

\begin{prop} \label{prop:conj}
Let $q$ be as above. Then $q$ divides~$\ansha$.
If
the Birch and Swinnerton-Dyer conjecture~(\ref{bsdform}) is true, then
$\ord_q \Big(\frac{\LAf}{\OAf} \Big)$ 
and~$\ord_q \big( \Mid \Sha(A_f) \miD \big)$ 
are both positive even numbers.
\end{prop}
\begin{proof}
By~\cite[Theorem~B]{emerton:optimal} (and considering that the order
of the cuspidal subgroup of~$J_0(N)$ is 
the numerator of~$\frac{N-1}{12}$ when $N$ is prime),
$q$ does not divide
$\prod_{\scriptscriptstyle{p |N}}  c_p(A_f)$ or
      ${ \Mid A_f(\Q) \miD \cdot \Mid A_f^{\vee}(\Q) \miD }$.
Thus if $q$ divides
$\frac{L(A_f,1)}{\OAf}$ then $q$ divides~$\Mid \Sha(E) \miD_{\rm an}$.
Now assume the Birch and Swinnerton-Dyer conjecture~(\ref{bsdform}),
so that $q$ divides~$\Mid \Sha(E) \miD$.
As mentioned towards the end of~\S7.3 in~\cite{dsw},
if $A^\vee_f[\qq]$ is irreducible for all maximal ideals~$\qq$
of~$\T$ with residue field of characteristic~$q$, 
then the $q$-primary part of~$\Sha(A^\vee_f)$
(and hence that of~$\Sha(A_f)$)
has order a perfect square. In our case, this irreducibility holds 
by~\cite[Prop.~14.2]{mazur:eisenstein}, and thus
$\ord_q \big( \Mid \Sha(A_f) \miD \big)$ 
is a positive even number.
Moreover, as mentioned above, $q$ does not divide any of the quantities
other than $\Mid \Sha(A_f) \miD$
on the right side of~(\ref{bsdform}), and hence we see that
$\ord_q \Big(\frac{\LAf}{\OAf} \Big)$ 
is a  positive even number.
\end{proof}

In particular, by Proposition~\ref{prop:conj} and its proof, 
if an odd prime~$q$ divides
$\frac{\LAf}{\OAf}$ or~$\ansha$, 
but does not divide the numerator of~$\frac{N-1}{12}$,
then $q^2$ (not just $q$) is expected
to divide $\frac{\LAf}{\OAf}$ and~$\ansha$.


Let $-D$ be a negative fundamental discriminant, and
as before, let $\epD = 
(\frac{-D}{\cdot})$ denote the associated quadratic character.
Suppose that $D$ is coprime to~$N$.
Then $f \tensor \epD$ is a modular form of level~$ND^2$. 
By a refinement of a theorem Waldspurger (see~\cite{luo-ram}), 
there exist infinitely many prime-to-$N$ discriminants~$-D$ such 
that $L(A_{f \tensor \epD},1) \neq 0$. Suppose $D$ is such that
$L(A_{f \tensor \epD}, 1) \neq 0$.
By Proposition~\ref{prop:ftensore}, the quantity
$\frac{L(A_{f \tensor \epD},1)}{\OAfmD}$ 
is an integer up to a power of~$2$, so it makes sense
to ask if an odd prime divides it.
Also, if $A_f$ is an elliptic curve and $(A_f, -D)$ satisfy
hypothesis~(**) mentioned at the beginning of Section~\ref{sec:twists}, 
then by Proposition~\ref{prop:twist},
$\frac{L(A_{f \tensor \epD},1)}{\OAfmD}$ is the 
algebraic part of the special $L$-value of the twist of~$A_f$ by~$-D$,
up to a power of~$2$.

\begin{thm} \label{sqthm}
Recall that the level~$N$ is assumed to be prime, and 
$q$ is an odd prime which does not divide the numerator of~$\frac{N-1}{12}$,
but divides $\frac{L(A_f,1)}{\OAf}$. Assume that $q$ satisfies
the following hypothesis:\\
(*) there exists 
a negative fundamental discriminant~$-D$ that is coprime to~$N$
such that $L(A_{f \tensor \epD},1) \neq 0$ and 
$q$ does not divide $\frac{L(A_{f \tensor \epD},1)}{\OAfmD}$.\\
Then $\ord_q \Big(\frac{\LAf}{\OAf} \Big)$ and~$\ord_q \big(\ansha \big)$ 
are both positive and even.
\end{thm}

We shall prove Theorem~\ref{sqthm} in Section~\ref{sec:proofs}. 
Assuming hypothesis~(*),
in view of Proposition~\ref{prop:conj}, Theorem~\ref{sqthm}
provides theoretical
evidence towards the Birch and Swinnerton-Dyer conjectural formula~(\ref{bsdform}).
We will say more about the hypothesis~(*) later in this section.

\begin{prop} \label{prop:sq}
Recall again that the level~$N$ is assumed to be prime. 
Suppose $q$ 
is an odd prime that does not divide the numerator of~$\frac{N-1}{12}$ and 
there is a normalized eigenform~$g \in S_2(\Gamma_0(N), \C)$ such that
$L(A_g,1) = 0$ and 
$f$ is congruent to~$g$ modulo a prime ideal over~$q$ 
in the ring of integers
of the number field generated by the Fourier coefficients
of~$f$ and~$g$. \\
(i) If the first part of the Birch and Swinnerton-Dyer conjecture 
is true for~$A_g$, then $q^2$ divides~$\Mid \Sha(A_f) \miD$.  \\
(ii) Suppose $q$ satisfies hypothesis~(*) of Theorem~\ref{sqthm}.
Then $q^2$ divides $\frac{L(A_f,1)}{\OAf}$ and 
the Birch and Swinnerton-Dyer conjectural value of~$\Mid \Sha(A_f) \miD$. 
In particular $\frac{L(A_f,1)}{\OAf} \equiv \frac{L(A_g,1)}{\OAg} \bmod q^2$.
\end{prop}

\begin{proof}
If the first part of the Birch and Swinnerton-Dyer conjecture (on rank)
is true for~$A_g$, then considering that
$L(A_g,1) = 0$, we see that $A_g$ has positive Mordell-Weil rank.
Part (i) now follows from~\cite[Thm~6.1]{agmer}.
By~\cite[Prop.~1.5]{agmer},
the hypotheses of
the proposition imply that $q$
divides $\LAf / \OAf$.
Thus part (ii) follows from the Theorem above.
\end{proof}

Subject to hypothesis~(*),
the proposition above shows some consistency between the predictions
of the two parts of the Birch and Swinnerton-Dyer conjecture.
There is a general philosophy that congruences
between eigenforms should lead to congruences between algebraic
parts of the corresponding special $L$-values, and there are 
theorems in this direction
(see~\cite{vatsal:canonical} and the references therein
for more instances). However, these theorems prove congruences
modulo primes, but not their powers. 
To our knowledge, part~(ii) of Proposition~\ref{prop:sq} above
is the first result of a form in which 
the algebraic parts of the special
$L$-value are congruent modulo 
the {\it square} of a congruence prime. 
In the rest of this section,
we give some heuristic and computational evidence for why 
hypothesis~(*) might always hold when $A_f$ is an elliptic curve, which
we denote by~$E$.
Suppose that $(E,-D)$ satisfies the hypothesis~(**) mentioned
at the beginning of Section~\ref{sec:twists}.
Then, by Proposition~\ref{prop:twist},
$\frac{L(A_{f \tensor \epD},1)}{\OAfmD}$
is the special $L$-value of the twisted elliptic curve~$E_{\scriptscriptstyle{-D}}$
up to a power of~$2$.

As mentioned before, 
by~\cite[Theorem~B]{emerton:optimal}, 
$q$ does not divide the orders of the arithmetic component
groups of~$E$, and hence by~\cite[Lem.~2.1]{prasanna:padic},
$q$ does not divide the orders of the arithmetic component
groups of~$E_{\scriptscriptstyle{-D}}$ either. Thus 
if one assumes the second
part of the Birch and Swinnerton-Dyer conjecture for~$E_{\scriptscriptstyle{-D}}$,
then the only way $q$ can divide $\frac{L(A_{f \tensor \epD},1)}{\OAfmD}$
is if $q$ divides the order of~$\Sha(E_{\scriptscriptstyle{-D}})$.

Now there is no clear reason for $q$ to divide the order 
of~$\Sha(E_{\scriptscriptstyle{-D}})$ for {\em every}~$D$.
Kolyvagin has asked whether for a given elliptic curve~$A$
and a prime~$q$, there is a twist of~$A$ such 
that $q$ does not divide the order of the Shafarevich-Tate group of the twist
(see Question~A in~\cite{prasanna:padic}). 
We are interested in the same
question, but with the added
restrictions that the level~$N$ is prime,
the special $L$-value
of the twist is nonzero,
and that $(E, -D)$ satisfies the hypothesis~(**). 

We now report on what numerical data we could gather regarding
this question.
Since we do not know a general algorithm to compute the actual order of the
Shafarevich-Tate group of an elliptic curve, we shall instead consider
the analytic orders and assume 
the second part of the Birch and Swinnerton-Dyer
conjecture to pass from analytic orders of the Shafarevich-Tate groups to their
actual orders.

Using the mathematical software sage, with its inbuilt Cremona's database 
for all elliptic curves of conductor up to~$130000$,
we considered all tuples $(N, E, p)$ such 
that $N$ is an integer less than $130000$, 
$E$ is an elliptic curve of conductor~$N$ 
with $\Mid \Sha(E) \miD_{\rm an}$ divisible by an odd prime, and
$p$ is an odd prime that divides $\Mid \Sha(E) \miD_{\rm an}$. 
For each such tuple, we looked for 
a negative fundamental discriminant  $-D$ such that 
$L(\ED, 1) \neq 0$, $ND^2 < 130000$ (to stay within  the range
of Cremona's database), and $D$
is coprime
to the discriminant of a chosen
Weierstrass equation 
$y^2 = x^3 + Ax + B$ of~$E$ with $A, B \in \Z$.
If we insisted on $N$ being prime, then
we found four tuples $(N, E, p)$ as above; for two of them,
we were able to find a~$D$ as above, in both of which
$p$ did not divide
$\Mid \Sha(\ED) \miD_{\rm an}$.
If we allow $N$ to be arbitrary, then
we found $357$ tuples $(N, E, p)$ as above, and for $103$ of them,
we were able to find a~$D$ as above, among which in $101$ cases,
$p$ did not divide
$\Mid \Sha(\ED) \miD_{\rm an}$. 
Of course, for the examples where we could not find a suitable~$D$
in the range of Cremona's tables, one may have to look beyond 
$ND^2 = 130000$ to satisfy hypothesis~(*). 
Indeed, even for $N$ as low as~$681$, which is the first level
at which an elliptic curve has the analytic order of the 
Shafarevich-Tate group divisible by an odd prime, 
the number of negative fundamental discriminants~$-D$ such that
${\rm gcd}(N, D) =1$ and $ND^2 < 130000$ is just~$4$.
In any case, when we could find a~$D$ satisfying the requirements above,
it was often the case that 
$p$ did not divide
$\Mid \Sha(\ED) \miD_{\rm an}$. 
Thus the data above does encourage the belief
that hypothesis~(*) always holds for elliptic curves
(even for non-prime levels). 
For more general newform quotients~$A_f$, we do not know how to do
computations (but see the remark at the end of Section~\ref{sec:overK}).

\comment{
As another test, we looked at all tuples $(N, E, p, D)$ satisfying
the above conditions. 
If we insisted on $N$ being prime, then
there were
only two tuples $(N, E, p, D)$ as above (they had different~$N$'s), and 
for both
$p$ did not divide
$\Mid \Sha(\ED) \miD_{\rm an}$. 
If we allow $N$ to be arbitrary, then we found
$103$ tuples $(N, E, p, D)$ as above
and in $101$ of them,
$p$ did not divide
$\Mid \Sha(\ED) \miD_{\rm an}$. 
Again, for the other two examples, where
$p$ did divide
$\Mid \Sha(\ED) \miD_{\rm an}$, we may have to go further beyond
the range of Cremona's database to find a suitable~$D$. 

We remark that Cremona has checked that the Manin constant
of all optimal elliptic curves of conductor up to~$130000$ is one
(see the appendix to~\cite{agst:manin}).
}

As mentioned above,
we have to assume 
the second part of the Birch and Swinnerton-Dyer
conjecture to pass from analytic orders of the Shafarevich-Tate groups to their
actual orders.
The careful reader would have noticed that we want to apply hypothesis~(*) to
give evidence for the second part of the Birch and Swinnerton-Dyer
conjecture, and at the same time we are assuming the conjecture to give
some credence to the hypothesis. While this may sound like circular reasoning,
the point is that the conjecture is being applied in different contexts, and
also our reasoning is not intended in any way to be a part of a proof.

One would of course hope that hypothesis~(*) is proved independent of 
the second part of the Birch and Swinnerton-Dyer conjecture.
While it is known that hypothesis~(*) does hold for 
all but finitely many primes~$q$ (e.g., see~\cite[Cor.~1]{ono-skinner:fourier}), it is not clear what that
finite list of primes is. 
Also, in~\cite[p.167-168]{bruinier-ono:coeff}, one finds a criterion for how 
big~$q$ needs to be, but the period they use (cf.~\cite[\S5]{bruinier:nonvanishing})
differs from the period we use by an unknown algebraic number 
(cf. the discussion in~\cite[Cor.~2]{kohnen:fourier}, 
and~\cite[Conj.~5.1]{prasanna:padic}). Thus unfortunately the theoretical results
mentioned in this paragraph do not help us much regarding hypothesis~(*).

\section{Special L-values over quadratic imaginary fields} \label{sec:overK}
Let $N$ be a  positive integer.
Let $f$ be a newform of weight~$2$ on~$\Gamma_0(N)$, and as before
let $A_f$ denote the associated newform quotient of~$J_0(N)$ 
over~$\Q$.
In this section, when $N$ is prime,
we give a formula for the algebraic part of the special
$L$-value of~$A_f$ over quadratic imaginary fields~$K$
in terms of the free abelian group on isomorphism classes of supersingular
elliptic curves in characteristic~$N$ (equivalently over conjugacy classes
of maximal orders in the definite quaternion algebra over~$\Q$ ramified at~$N$
and~$\infty$)
which shows that this algebraic part is a 
perfect square away from the prime two and 
the primes dividing the discriminant
of~$K$.

We start by recalling the definition of the ``archemedian volume'' $\Omega(A_f/K)$
alluded to in the introduction.
Let $d = \dim A_f$ and  
let $F$ be a number field. 
Let
$\omega_1, \ldots, \omega_d$ be a basis of~$H^0(A_f,\Omega^1_{A_f/\Q})$
associated to a $\Z$-basis of~$S_2(\Gamma_0(N), \Z)[I_f]$.
Then $\omega_1, \ldots, \omega_d$ 
is also an $F$-basis of~$H^0(A_f,\Omega^1_{A_f/F})$.
Let $W$ denote the group of invariant differentials
on the N\'eron model~$\NerA_{\OO}$
of~$A_f$ over~$\OO$, the ring of integers of~$F$.
Then $\wedge^d W =  \cF \cdot \wedge_i \omega_i$ 
for some fractional ideal~$\cF$ of~$\OO$
(cf.~\cite[\S~III.5]{lang:nt3}).
We will call the ideal~$\cF$ the {\em Manin ideal} of~$A_f$ over~$F$.
If $F = \Q$, then the absolute value of a generator
of the Manin ideal is just the Manin constant of~$A_f$
(as defined in~\cite{agst:manin})
and is denoted~$\cAf$. If $A_f$ is an elliptic curve,
then this definition of the Manin constant agrees with 
the one given in Section~\ref{sec:twists}
for optimal elliptic curves. 
The Manin constant~$\cAf$ is conjectured
to be one; 
it is known that $\cAf$ is an integer, and 
if $p$ is a prime such that $p^2 \nmid 2N$, then 
$p$ does not divide $\cAf$ (see~\cite{agst:manin} for details).

\begin{lem} \label{lem:manin}
The Manin ideal~$\cF$ is supported on the set of maximal ideals~$\m$
of~$\OO$ such that the residue characteristic of~$\m$ divides
either~$\cAf$ or the discriminant of~$\OO$.
\end{lem}
\begin{proof}
Suppose $\m$ is a maximal ideal
of~$\OO$ such that the residue characteristic~$\ell$ of~$\m$ divides
neither~$\cAf$ nor the discriminant of~$\OO$.
By~\cite[\S7.2, Cor.~2]{neronmodels}, over discrete valuation rings,
the formation of N\'eron models is compatible with unramified extensions.
Thus, considering that
$\ell$ is coprime to the discriminant of~$\OO$,
 $H^0(\NerA_{\OO}, \Omega_{\NerA_{\OO}/\OO}) \tensor_\OO \OO_\m 
= H^0(\NerA_{\OO_\m}, \Omega_{\NerA_{\OO_\m}/{\OO_\m}}) 
= H^0(\NerA_{\Z_\ell}, \Omega_{\NerA_{\Z_\ell}/{\Z_\ell}}) \tensor_{\Z_\ell} \OO_\m
= H^0(\NerA_{\Z}, \Omega_{\NerA_{\Z}/{\Z}}) \tensor_{\Z} 
\Z_\ell \tensor_{\Z_\ell} \OO_\m 
= H^0(\NerA_{\Z}, \Omega_{\NerA_{\Z}/{\Z}}) \tensor_{\Z} \OO_\m
$. Thus, a $\Z$-basis 
$\omega'_1, \ldots, \omega'_d$ of~$H^0(\NerA_{\Z}, \Omega_{\NerA_{\Z}/{\Z}})$
is a $\OO_\m$-basis 
of~$H^0(\NerA_{\OO}, \Omega_{\NerA_{\OO}/\OO}) \tensor_\OO \OO_\m$.
Since $\ell$ does not divide~$\cAf$,
we see that $(\wedge^d \omega'_i) \tensor \Z_\ell 
=  (\wedge^d \omega_i) \tensor \Z_\ell $. 
Hence $(\wedge^d \omega'_i) \tensor \OO_\m
=  (\wedge^d \omega_i) \tensor \OO_\m$. 
In view of all this, it follows that
$\cF \tensor_\OO \OO_\m$ is trivial, and the lemma follows.
\end{proof}
\comment{Let
$\omega_1, \ldots, \omega_d$ be a basis of~$H^0(A_f,\Omega^1_{A_f/\Q})$;
it is also a $K$-basis of~$H^0(A_f,\Omega^1_{A_f/K})$.
Let $W$ denote the group of invariant differentials
on the N\'eron model~$\NerA_{\OO_K}$
of~$A_f$ over~$\OO_K$, the ring of integers of~$K$.
Then $\wedge^d W = \cF \wedge_i \omega_i$ 
for some fractional ideal~$\cF$ of~$\OO_K$.
}
Let $c_1, \ldots, c_{2d}$ be a basis of~$H_1(A_f(\C),\Z)$.
The {\em complex period matrix} of~$A_f$ (with respect to
the chosen basis) is the $2d \times 2d$ matrix
$A = (\int_{c_i} \omega_j, \int_{c_i} \overline{\omega}_j)$.
Recall that  $K$ is a quadratic imaginary field; let $-D$
be its discriminant.
The ``archimedean volume'' of~$A_f$ over~$K$ is
\begin{eqnarray} \label{eqn:comper}
\Omega(A_f/K) = | \det (A) | \cdot N^K_\Q(\cF) / D^{d/2}
\end{eqnarray}
(this coincides with the definition of~$C_{A,\infty}$ 
in~\cite[\S~III.5]{lang:nt3}).
\comment{
Now in view of the fact that the 
Manin constant, denoted $\cF$ in loc. cit., is a power of 2
by~\cite[Cor.~4.1]{mazur:rational}, 
considering that the level~$N$ is prime).
}

Let $N$ be prime in the rest of this section.
We next give a formula for the ratio $\frac{L(A_f/K,1)}{\Omega(A_f/K)}$,
which is the left hand side of the 
Birch and Swinnerton-Dyer
conjectural formula (Conjecture~\ref{bsd}) for~$A_f$ over~$K$. 

Let $\{E_0, E_1, \ldots, E_g\}$ be a set of representatives
for the isomorphism classes of supersingular elliptic
curves in characteristic~$N$, where $g$ is the genus of~$X_0(N)$.
We denote the class of~$E_i$ by~$[E_i]$.
Let $\sP$ denote the divisor group supported on 
the $[E_i]$ and let $\Pz$ denote the subgroup
of divisors of degree~$0$.
For $i = 1, 2, \ldots, g$, 
let $R_i = {\rm End\ } E_i$. Each $R_i$ is a maximal order in the definite
quaternion algebra ramified at~$N$ and~$\infty$,
which we denote by~$\mathcal{B}$ and in fact, the $R_i$'s are
representatives of the conjugacy
classes  of maximal orders of~$\mathcal{B}$.
Moreover, setting $I_i = {\rm Hom}(E_0,E_i)$, we see that 
the $I_i$ are representatives for the isomorphism classes of right
locally free rank one right modules over~$R_0$. 
Let $\OD$ denote the quadratic order of
discriminant~$-D$, $h(-D)$ the number of classes
of~$\OD$, $u(-D)$ the order of $\OO^*_{-D}/\langle \pm 1 \rangle$,
and $h_i(-D)$ the number of optimal embeddings of~$\OD$
in~$R_i$ modulo conjugation by~$R_i^*$.
Following~\cite{gross:heights}, we define
$$\chi_{\scriptscriptstyle{D}}
 = \frac{1}{2 u(-D)} \sum_{i=0}^g h_i(-D)[E_i] \in \sP \tensor \Q.$$
Let $w_i = \Mid {\rm Aut} E_i \miD = \Mid R_i^*/\langle \pm 1 \rangle \miD$. 
Define the Eisenstein element in~$\sP \tensor \Q$ as
$a_{\scriptscriptstyle{E}} = \sum_{i=0}^g \frac{[E_i]}{w_i}$.
Let $\chiD = \chi_{\scriptscriptstyle{D}} - 
\frac{12}{p-1} \deg(e_{\scriptscriptstyle{D}}) 
a_{\scriptscriptstyle{E}}$. 
Let \mbox{$n = {\rm numr}(\frac{N-1}{12})$}; then
$n \chiD \in \Pz$. 
Since the level~$N$ is prime, the Hecke algebra~$\T$ is semi-simple,
and hence we have an isomorphism $\T \tensor \Q \isom \T/I_f \tensor \Q
\oplus B$ of $\T \tensor \Q$-modules for some
$\T \tensor \Q$-module~$B$.
Let $\pi$ denote element of $\T \tensor \Q$ that is
the projection on the first factor.
We prove the following in Section~\ref{sec:proofs}:

\begin{thm} \label{prop:sqoverk}
Recall that the level~$N$ is prime. Let $K$ be a quadratic imaginary field
with discriminant~$-D$  that is
coprime to~$N$.
If $L(A_f/K,1) \neq 0$, then up to powers of primes dividing~$2D$,
$$\frac{L(A_f/K,1)}{\Omega(A_f/K)} = 
\frac{\mid \pi ({\Pz}) : \pi({\T n \chiD})  \mid^2}
{N^K_\Q(\cK) \cdot n^2}\ \ \ .$$ 
Moreover, $\frac{L(A_f/K,1)}{\Omega(A_f/K)}$ is a perfect square
up to powers of primes dividing~$2D$. 
\end{thm}

This addresses the issue raised in~\cite[p.~236]{marusia:module} that
as of the writing of loc. cit., one did not have a way of expressing special
$L$-values over~$K$
in terms of the module~$\sP$. Also, it may be possible to use
the formula above for computations using Brandt matrices (cf.~\cite{kohel:hecke}). 
Note that up to powers of primes dividing~$2D$,
$\frac{L(A_f/K,1)}{\Omega(A_f/K)}$ equals
$\frac{\LAf}{\OAf} \cdot \frac{L(A_{f \tensor \epD},1)}
{\OAfmD}$
(see formula~(\ref{eqn:prod}) in Section~\ref{sec:proofs}). 
Thus if the formula in Theorem~\ref{prop:sqoverk} could be used
for computations, then considering that one already knows how to compute
$\frac{\LAf}{\OAf}$ (see~\cite[\S4]{agst:bsd}), one could compute 
$\frac{L(A_{f \tensor \epD},1)}{\OAfmD}$ systematically and check whether
the hypothesis~(*) 
of Theorem~\ref{sqthm}
holds in particular examples for odd primes~$q$ not dividing~$D$.

\section{Proofs of Theorems~\ref{sqthm} and~\ref{prop:sqoverk}} 
\label{sec:proofs}

In this section, we prove Theorems~\ref{sqthm} and~\ref{prop:sqoverk}.
We shall be using results 
from~\cite{marusia:module}, and details of some
of the facts that we use here routinely may be found in loc. cit.

Let $\mathcal{H}$ denote the complex upper half plane, and
let ${\{0,i\infty\}}$ denote the projection of the geodesic
path from~$0$ to~$i\infty$ in~${\mathcal{H}} \cup {\PP}^1 ({\Q})$
to~$X_0(N)({\C})$.
We have an isomorphism
$$H_1(X_0(N),{\Z}) \tensor {\R} \stackrel{\cong}{\lra}
{\rm Hom}_{{\C}} (H^0(X_0(N), \Omega^1),{\C}),$$ obtained
by integrating differentials along cycles. 
Let $e$ be the element of $H_1(X_0(N),{\Z}) \tensor {\R}$ that corresponds
to the map $\omega \mapsto - \int_{\{0,i\infty\}} \omega$ under this
isomorphism. It is called the {\em winding element}.
By the Manin-Drinfeld Theorem,
(see~\cite[Chap. {I}{V}, Theorem~$2.1$]{lang:modular} 
and~\cite{manin:parabolic}),
\mbox{$e \in H_1(X_0(N),\Z) \tensor \Q$}. Also, since the complex
conjugation involution on $H_1(X_0(N),\Z)$ is induced by the map
$z \mapsto - \overline{z}$ on the complex upper half plane,
we see that $e$ is invariant under complex conjugation.
Thus $e \in H_1(X_0(N),\Z)^+ \tensor \Q$.
Let $H^+$ and~$H^-$ denote the subgroup of elements of~$H_1(X_0(N),\Z)$
on which the complex conjugation
involution acts as~$1$ and~$-1$ respectively.  

Assume henceforth that $N$ is prime (which is a hypothesis for
the theorems that we want to prove).
Consider the $\T[1/2]$-equivariant isomorphism 
\begin{eqnarray} \label{phiisom}
\Phi: \Pz[1/2] \tensor_{\T[1/2]} \Pz[1/2] \ra
H^+[1/2] \tensor_{\T[1/2]} H^-[1/2]
\end{eqnarray}
obtained from~\cite[Prop.~4.6]{marusia:module}
(which says that both sides of~(\ref{phiisom}) are isomorphic 
to~$S_2(\Gamma_0(N),\Z)[1/2]$, and whose proof
relies on results of~\cite{emerton:ss}).
By~\cite[Thm~0.2]{marusia:module}, we have
$\Phi_\Q(\chiD \tentq \chiD) = e \tentq \eD$,
where the subscript~$\Q$ stands for tensoring with~$\Q$
(this follows essentially from~\cite[Cor~11.6]{gross:heights},
along with its generalization~\cite[Thm~1.3.2]{zhang:gz}).
Thus $\Phi_\Q$ induces an isomorphism 
\begin{eqnarray} \label{eDtoetD}
\T[1/2] (n \chiD \tenth n \chiD)
\isom \T[1/2] n e \tenth \T[1/2] n \eD.
\end{eqnarray}
Note that $ne \in H^+$ by II.18.6 and II.9.7 of~\cite{mazur:eisenstein}.


Recall that since the level~$N$ is prime, 
the Hecke algebra~$\T$ is semi-simple,
and hence we have an isomorphism $\T \tensor \Q \isom \T/I_f \tensor \Q
\oplus B$ of $\T \tensor \Q$-modules for some
$\T \tensor \Q$-module~$B$.
Recall also that  $\pi$ denotes the element of $\T \tensor \Q$ that is
the projection on the first factor. 
In this section, if $X$ and~$Y$ are $\T$-modules with $Y \subseteq X$,
then we shall write $\big| \pi\big(\frac{X}{Y} \big) \big|$
for $\mid \pi(X) : \pi(Y) \mid$, which is an integer;
we are doing this so that the formulas do not look too terrible.

\begin{prop} \label{hpart}
\begin{eqnarray*}
\bigg| \pi \bigg( \frac{H^+[1/2]}{\T[1/2] n e} \bigg) \bigg|
\cdot \bigg| \pi \bigg( \frac{H^-[1/2]}{\T[1/2] n \eD} \bigg) \bigg|
= \bigg| \pi \bigg( \frac{H^+[1/2] \tent H^-[1/2])}{\T[1/2] n e \tent 
\T[1/2] n \eD} \bigg) \bigg|
\end{eqnarray*}
\end{prop}
\begin{proof}
By~\cite[\S15]{mazur:eisenstein}, if $\mm$ is a Gorenstein
maximal ideal of~$\T$ with odd residue characteristic, 
then $H^+_\mm$ and~$H^-_\mm$ are free $\T_\mm$-modules of
of rank one. 
Since the level is prime, the only non-Gorenstein ideals
are the ones lying over~$2$, a prime that we are systematically
inverting anyway.

Let $\mm$ be a maximal ideal of~$\T$ with odd residue characteristic. 
Let $x$ be a generator of $\Hpm$ as a free $\T_\mm$-module,
and let $y$ be a generator of $\Hmm$ as a free $\T_\mm$-module.
Then there exists $t_1 \in \Tm$ such that $n e = t_1 x$ 
and $t_2 \in \Tm$ such that $n \eD = t_2 y$. 
We have
\comment{
\noindent $|\pi (\Hpm \tentm \Hmm)/ 
  \pi(\Tm n e \tentm \Tm n \eD)| 
= |\pi (\Tm x \tentm \Tm y)/ 
  \pi(\Tm t_1 x \tentm \Tm t_2 y)| \\
= |\pi (\Tm (x \tentm y))/ 
  t_1 t_2 \pi(\Tm (x \tentm y))|
= |\pi (\Tm)/ t_2 t_1 \pi(\Tm)|\\
= |\pi (\Tm)/ t_1 \pi(\Tm)|\cdot |\pi (t_1 \Tm)/ t_2 \pi(t_1 \Tm)|.
$ \\
}
\begin{eqnarray*}
\bigg| \frac{\pi (\Hpm \tentm \Hmm)}{   \pi(\Tm n e \tentm \Tm n \eD)} \bigg| 
& = & 
\bigg| \frac{\pi (\Tm x \tentm \Tm y)}{  \pi(\Tm t_1 x \tentm \Tm t_2 y)} \bigg| \\
& = & \bigg| \frac{\pi (\Tm (x \tentm y))}{   t_1 t_2 \pi(\Tm (x \tentm y))} \bigg|\\
& = & \bigg| \frac{\pi (\Tm)}{ t_2 t_1 \pi(\Tm)} \bigg| \\
& = &\bigg| \frac{\pi (\Tm)}{ t_1 \pi(\Tm)}  \bigg|
\cdot \bigg| \frac{\pi (t_1 \Tm)}{ t_2 \pi(t_1 \Tm)} \bigg|.
\end{eqnarray*}

\noindent{\em Claim:} 
$$\bigg| \frac{\pi (t_1 \Tm)}{ t_2 \pi(t_1 \Tm)} \bigg| = \bigg| \frac{\pi (\Tm)}{ t_2 \pi(\Tm)} \bigg|\ \ .$$
\begin{proof}
Consider the map 
$\psi: \pi (\Tm) \ra \pi (t_1 \Tm)/t_2 \pi(t_1 \Tm)$
given as follows: if $t \in \Tm$, then $\pi(t) \mapsto \pi(t_1 t)$.
If $\pi(t)$ is in the kernel of~$\psi$, then
$\pi(t_1 t) = \pi(t_2 t_1 t')$ for some $t' \in \Tm$. Then
$\pi(t_1(t - t_2 t')) = 0$, and since $\pi(t_1) \neq 0$
(as $L(A_f,1) \neq 0$), we have $\pi(t) = \pi(t_2 t')$.
Thus the kernel of~$\psi$ is $t_2 \pi(\Tm)$, which proves
the lemma.
\end{proof}

Using the claim and the series of equalities above, we have\\
\begin{eqnarray*}
\bigg| \frac{\pi (\Hpm \tentm \Hmm)}{ 
  \pi(\Tm n e \tentm \Tm n \eD)} \bigg| 
& = & \bigg| \frac{\pi (\Tm)}{ t_1 \pi(\Tm)} \bigg| \cdot \bigg| \frac{\pi (\Tm)}{ t_2 \pi(\Tm)} \bigg|\\
& = & \bigg| \frac{\pi (\Tm x)}{ t_1 \pi(\Tm x)} \bigg| \cdot \bigg| \frac{\pi (\Tm y)}{ t_2 \pi(\Tm y)} \bigg|\\
& = & \bigg| \frac{\pi(\Hpm)}{\pi(\Tm n e)} \bigg| \cdot \bigg| \frac{\pi(\Hmm)}{ \pi(\Tm n \eD)} \bigg|\\
& = & \bigg| \pi \bigg( \frac{\Hpm}{\Tm n e} \bigg) \bigg|
\cdot \bigg| \pi \bigg( \frac{\Hmm}{\Tm n \eD} \bigg) \bigg| \ \ .
\end{eqnarray*}

Since this is true for every~$\mm$ with odd residue characteristic,
we get the statement in the proposition.
\end{proof}

\begin{prop} \label{ppart}
$$
 \bigg| \pi \bigg( \frac {\Pz[1/2] \tenth 
\Pz[1/2]} {\T[1/2] (n \chiD \tenth n \chiD)} \bigg) \bigg|
= \bigg| \pi \bigg (\frac {\Pz[1/2]}{\T[1/2] n \chiD} \bigg) \bigg|^2.
$$
\end{prop}
\begin{proof}
By~\cite[Thm~0.5]{emerton:ss}, if $\mm$ is a Gorenstein
maximal ideal of~$\T$, then $\Pzm$ is a free $\T_\mm$-module
of rank one; let $x$ be a generator. Then $n \chiD = t x$ for some $t \in \Tm$.
Hence in a manner similar to the steps in the proof of
Proposition~\ref{hpart}, we have
$$ \bigg| \pi \bigg( \frac {\Pzm \tenth \Pzm} 
        {\Tm (n \chiD \tenth n \chiD)} \bigg) \bigg|
= \bigg| \pi \bigg( \frac {\Tm x \tenth \Tm x} 
        {\Tm (t x \tenth tx)} \bigg) \bigg|
= \bigg| \pi \bigg (\frac {\Tm}{t^2 \Tm} \bigg) \bigg|$$
$$= \bigg| \pi \bigg (\frac {\Tm}{t \Tm} \bigg) \bigg|^2
= \bigg| \pi \bigg (\frac {\Pzm}{\Tm n \chiD} \bigg) \bigg|^2. $$
Since this holds for every maximal ideal~$\mm$ of odd residue
characteristic, we get the proposition.
\end{proof}

By formula~(\ref{phiisom}), formula~(\ref{eDtoetD}), 
Proposition~\ref{hpart}, and Proposition~\ref{ppart}, we have
\begin{eqnarray} \label{semifinaleqn}
\bigg| \pi \bigg( \frac{H^+[1/2]}{\T[1/2] n e} \bigg) \bigg|
\cdot \bigg| \pi \bigg( \frac{H^-[1/2]}{\T[1/2] n \eD} \bigg) \bigg| 
= \bigg| \pi \bigg (\frac {\Pz[1/2]}{\T[1/2] n \chiD} \bigg) \bigg|^2
\end{eqnarray}

Let $\OAfp = {\rm disc}(H_1(A_f,\Z)^+ \times S_f \ra {\C})$;
it differs from~$\Omega(A_f)$ by a power of~$2$
(by~\cite[Lemma~2.4]{agmer}).
By the proof of Theorem~2.1 of~\cite{agmer}, we have
$$\bigg| \pi \bigg( \frac{H^+}{\T (n e)} \bigg) \bigg|
= 
n \cdot \frac{\LAf}{\OAfp}.$$
Using this and Proposition~\ref{prop:ftensore}, equation
(\ref{semifinaleqn}) says that up to a power of~$2$, 
\begin{eqnarray} \label{finaleqn}
\frac{\LAf}{\OAfp} 
\cdot \frac{L(A_{f \tensor \epD},1)}{\OAfmD}
= \frac{1}{n^2} \cdot 
\bigg| \pi \bigg (\frac {\Pz[1/2]}{\T[1/2] n \chiD} \bigg) \bigg|^2.
\end{eqnarray}

\begin{proof}[Proof of Theorem~\ref{prop:sqoverk}]
We have  $L(A_f/K, s) = L(A_f,s) \cdot L(A_{f \tensor \epD},s)$,
and 
by Corollary~\ref{cor:periods} in Section~\ref{sec:app},
we have $\Omega(A_f/K) = N^K_\Q(\cK) \cdot \OAfp \cdot \OAfm / (-D)^{d/2} $,
up to a sign. 
Thus we have
\begin{eqnarray} \label{eqn:prod}
\frac{L(A_f/K,1)}{\Omega(A_f/K)} = 
\frac{1}{N^K_\Q(\cK)} \cdot \frac{\LAf}{\OAfp} \cdot \frac{L(A_{f \tensor \epD},1)}{\OAfmD/(-D)^{d/2}} \ ,
\end{eqnarray}
up to a sign.
The first claim of
Theorem~\ref{prop:sqoverk} now follows from~(\ref{finaleqn}).
The second claim follows from the first since
$N^K_\Q(\cK)$ is coprime to $2D$ by Lemma~\ref{lem:manin}, considering
that by~\cite[Cor.~4.1]{mazur:rational}, $\cAf$ is a power of~$2$ 
since $N$ is prime.

\end{proof}

\begin{proof} [Proof of Theorem~\ref{sqthm}]
If an odd prime~$q$ divides $\frac{\LAf}{\OAf}$
(which recall differs from $\frac{\LAf}{\OAfp}$ by a power of~$2$)
and $q$ does not divide 
$\frac{L(A_{f \tensor \epD},1)}{\OAfmD}$,
then by~(\ref{finaleqn}), $\ord_q \Big(\frac{\LAf}{\OAf} \Big)$
is even (and positive).
By~\cite[Theorem~B]{emerton:optimal},
we have $\Mid A_f(\Q) \miD = \Mid A_f^{\vee}(\Q) \miD$ and
this order divides the numerator of~$\frac{N-1}{12}$.
Thus if $q$ does not divide the numerator of~$\frac{N-1}{12}$,
then from~(\ref{bsdform}), $\ord_q \big(\ansha \big)$ is positive and even.
This proves Theorem~\ref{sqthm}.
\end{proof}

\section{Appendix: period matrices} \label{sec:app}

In this section, we give a formula for the determinant of 
the ``complex period matrix'' of an abelian variety. 
The result is probably well known, but we could not find
a suitable reference.

Let $Y$ be an abelian variety over~$\Q$ of dimension~$d$. 
Let $\omega_1, \ldots, \omega_d$ be a basis of~$H^0(Y,\Omega^1_{Y/\Q})$.
Let $c_1, \ldots, c_{2d}$ be a basis of~$H_1(Y(\C),\Z)$.
We define the associated {\em complex period matrix} of~$Y$
as the $2d \times 2d$ matrix
$A = (\int_{c_i} \omega_j, \int_{c_i} \overline{\omega}_j)$;
this matrix depends on the choices of the bases made above.

We have an action of complex conjugation~$c$ on~$Y(\C)$, and 
hence on $H_1(Y(\C),\Z)$.
Let $H_1(Y(\C),\Z)^+$ denote the subgroup of elements of~$H_1(Y(\C),\Z)$
that are fixed by complex conjugation, and let
$H_1(Y(\C),\Z)^-$ denote the subgroup of elements~$x$ of~$H_1(Y(\C),\Z)$
such that $c(x) = -x$. 
Let $\gamma_1, \ldots, \gamma_d$ be a basis of~$H_1(Y(\C),\Z)^+$,
and let $\gamma'_1, \ldots, \gamma'_d$ be a basis of~$H_1(Y(\C),\Z)^-$.
Let $B$ denote the $d \times d$
matrix whose $(i,j)$-th entry is $\int_{\gamma_i} \omega_j$
and let 
$C$ denote the $d \times d$
matrix whose $(i,j)$-th entry is $\int_{\gamma'_i} \omega_j$.

\begin{lem} \label{lem:periods} We have
$\det A = \det B \cdot \det C$ up to a sign and up to a power of~$2$.
\end{lem}
\begin{proof}
Let 
$A_{1,2}$ denote the $d \times d$ matrix whose
$(i,j)$-th entry is $\int_{\gamma_i} \overline{\omega}_j$, 
and 
let $A_{2,2}$ denote the $d \times d$ matrix whose
$(i,j)$-th entry is $\int_{\gamma'_i} \overline{\omega}_j$.
Consider the 
$2d \times 2d$ matrix 
$$A' = 
\left[
\begin{matrix}
B &  A_{1,2} \\
C & A_{2,2} 
\end{matrix}
\right]
$$
Now the set $\{ \gamma_1, \ldots, \gamma_d,
\gamma'_1, \ldots, \gamma'_d \}$ generates a subgroup 
of~$H_1(A(\C),\Z)$ of index a power of~$2$, and 
thus up to a sign and up to a power of~$2$, 
we have 
\begin{eqnarray} \label{eqn:det}
\det(A)  =  \det(A')\ . 
\end{eqnarray}
Now if $\gamma \in H_1(A(\C),\Z)$, and $\overline{\gamma}$ denotes its complex
conjugate, then for $j = 1, \ldots, d$, since $\omega_j$ is $\Q$-rational,
we have
$\int_\gamma \overline{\omega}_j = \int_{\overline{\gamma}} \omega_j$.
In particular, if $\gamma \in H_1(A(\C),\Z)^+$,
then 
$\int_\gamma \overline{\omega}_j = \int_{\gamma} \omega_j$,
and if $\gamma \in H_1(A(\C),\Z)^-$,
then 
$\int_\gamma \overline{\omega}_j = - \int_{\gamma} \omega_j$.
Thus we see that $A_{1,2} = B$ and $A_{2,2} = - C$.
Thus 
$$A' = 
\left[
\begin{matrix}
B &  B \\
C & -C
\end{matrix}
\right] \ .
$$
From this, we see that $\det(A') = -2 \det(B) \det(C)$.
The lemma now follows from~(\ref{eqn:det}).
\end{proof}

We remark that the discussion above holds
even if we replace $\Q$ by~$\R$ throughout.

\begin{cor} \label{cor:periods}
Let $N$ be a  positive integer.
Let $f$ be a newform of weight~$2$ on~$\Gamma_0(N)$, and as before
let $A_f$ denote the associated newform quotient of~$J_0(N)$ 
over~$\Q$.
Recall that $\OAfp = {\rm disc}(H_1(A_f,\Z)^+ \times S_f \ra {\C})$,
and $\OAfm = {\rm disc}(H_1(A_f,\Z)^- \times S_f \ra {\C})$,
where $S_f = S_2(\Gamma_0(N),\Z)[I_f]$. Let $K$ be a quadratic
imaginary field of discriminant~$-D$, and let
 $\Omega(A_f/K)$ be the ``complex period'' of~$A_f$ over~$K$ as defined
in formula~(\ref{eqn:comper}) of Section~\ref{sec:overK}. Then up to 
a sign,
$$\Omega(A_f/K) = N^K_\Q(\cK) \cdot \OAfp \cdot  \OAfm / (-D)^{d/2},$$
where $\cK$ is the Manin ideal of~$A_f$ over~$F$,
as defined at the beginning of Section~\ref{sec:overK}.
\end{cor}
\begin{proof}
Take $Y = A_f$ in the discussion at the beginning of this section, 
and take $\omega_1, \ldots, \omega_d$ to be the differentials 
in~$H^0(A_f,\Omega^1_{A_f/\Q})$ corresponding to a basis of~$S_f$.
Let the matrices
$A$, $B$, and~$C$ be as above, for the choices made
in the previous sentence.
Then by definition, 
$\Omega(A_f/K) = | \det (A) | / (-D)^{d/2}$, 
$\OAfp = \det(B)$, and 
$\OAfm = \det(C)$.
 
Now if $\gamma \in H_1(A(\C),\Z)$, and $\overline{\gamma}$ denotes its complex
conjugate, then for $j = 1, \ldots, d$, since $\omega_j$ is $\Q$-rational,
we have
$\overline{ \int_\gamma \omega_j} = \int_{\overline{\gamma}} \omega_j$.
In particular, if $\gamma \in H_1(A(\C),\Z)^+$,
then 
$\overline{ \int_\gamma \omega_j} = \int_{\gamma} \omega_j$,
so $\int_{\gamma} \omega_j$ is real. Hence all the entries of the matrix~$B$
are real. Hence $|\det(B)| = \det(B)$ up to a sign.
Similarly, if $\gamma \in H_1(A(\C),\Z)^-$,
then 
$\overline{ \int_\gamma \omega_j} =  - \int_{\gamma} \omega_j$,
so $\int_{\gamma} \omega_j$ is purely imaginary. 
Thus all the entries of the matrix~$C$ are purely imaginary. 
Hence $|\det(C)| = (\sqrt{-1})^d \det(C)$ up to a sign.

Thus by Lemma~\ref{lem:periods} and the discussion
in the two paragraphs above, we see
that  up to a sign,
$\Omega(A_f/K) = | \det (A) | / D^{d/2}
= |\det(B)| \cdot |\det(C)| / D^{d/2}
= \det(B)  \cdot \det(C)  \cdot (\sqrt{-1})^d / D^{d/2}
= \det(B)  \cdot  \det(C) / (-D)^{d/2}
= \OAfp \cdot  \OAfm / (-D)^{d/2}$, as was to be shown.
\end{proof}

\section*{Acknowledgments}

This paper owes its existence to 
Loic Merel. He pointed out to us that Gross' formula~\cite{gross:heights}
should have some implications for the squareness of the order of
the Shafarevich-Tate group, and also indicated the relevance
of~\cite{marusia:module}. 
The author's task was to work out the details
and see what precise implications
could be drawn. We would like to thank Loic Merel for suggesting 
this project as well as for several discussions regarding it.
We are also grateful to J.~Brown, W.~Li, K.~Ono, and K.~Prasanna for some
discussions.
Some of the computations were done on sage.math.washington.edu, which
is supported by National Science Foundation Grant No. DMS-0821725.

\bibliographystyle{amsalpha}         

\providecommand{\bysame}{\leavevmode\hbox to3em{\hrulefill}\thinspace}
\providecommand{\MR}{\relax\ifhmode\unskip\space\fi MR }
\providecommand{\MRhref}[2]{%
  \href{http://www.ams.org/mathscinet-getitem?mr=#1}{#2}
}
\providecommand{\href}[2]{#2}

\end{document}